\documentclass[11pt]{article}
\usepackage{amsmath,amsthm,amsfonts}
\theoremstyle{definition}

\theoremstyle{definition}

\theoremstyle{definition}
 
\begin{document}
\title{\LARGE\bf 
Triality, Structurable and Pre-structurable Algebras\\}
\date{}
\author{\Large Noriaki Kamiya$^1$ and Susumu Okubo$^2$\\ \\
$^1$Department of Mathematics, University of Aizu\\
965-8580,  Aizuwakamatsu, Japan\\ 
E-mail: kamiya@u-aizu.ac.jp\\ \\
$^2$Department of Physics and Astronomy, University of Rochester\\
Rochester, New York, 14627, U.S.A\\
E-mail: okubo@pas.rochester.edu}
\maketitle
\thispagestyle{empty}
\begin{abstract}
We introduce a notion of Pre-structurable Algebras based
upon  triality relations and study its relation to
structurable algebra of Allison, as well as to Lie algebras satisfying triality.
\end{abstract}
\par
\vskip 3mm
{\bf I Introduction, preliminary and  pre-structurable algebra}
\par
The structurable algebra [A], 
which is a class of nonassociative algebras has many interesting properties.
First,
it satisfies a triality relation [A-F].
Moreover,
we can construct a Lie algebra $L$
(see [A-F])
that any simple classical Lie algebra
can be constructed from some appropriate 
structurable algebra.
Further any such Lie algebra $L$ is invariant
under the symmetric group
$S_{4}$ and is a $BC_{1}$-graded Lie algebra of 
type $B_{1}$
(see [E-O.2] and [E-K-O]).
\par
Although the structurable algebra has been originally defined in term
of Kantor's triple system [Kan],
or equivalently of $(-1,1)$
Freudenhthal-Kantor triple system
[Y-O],
it can be also defined without any reference upon the
triple system [A-F].
Here,
following [O.1]
we base its definition in terms of
triality relations as follows:
\par
Let $(A,-)$ be an algebra over a field $F$ with bi-linear
product denoted by
juxtaposition  $ab\in A$ for
$a,b\in A$
with the unit element
$e$ and with the
involution map
$a\rightarrow {\bar a}$ and  $ \overline{ab}=\bar{b}\bar{a}$.
\par
Let
$$t_{j}:A\otimes A\rightarrow
\ {\rm End}\ A,
(j=0,1,2)\eqno(1.1)$$
be given by (see [A-F])
$$
t_{1}(a,b):=
l({\bar b})l(a)-
l({\bar a})l(b)\eqno(1.2a)$$
$$
t_{2}(a,b):=
r({\bar b})r(a)-
r({\bar a})r(b)\eqno(1.2b)$$
$$
t_{0}(a,b):=
r({\bar a}b-{\bar b}a)+
l(b)l({\bar a})-
l(a)l({\bar b})\eqno(1.2c)$$
where
$l(a)$ and
$r(a)$
are standard multiplication operaters given by
$$
l(a)b=ab\eqno(1.3a)$$
$$
r(a)b=ba.
\eqno(1.3b)$$
We then  see that $t_{0}(a,b)$ satisfies automatically
$$
t_{0}(a,b)c+
t_{0}(b,c)a+
t_{0}(c,a)b=0.\eqno(1.4)$$
Suppose now that
$t_{j}(a,b)$
satisfy the triality relation
$$
\overline{ t_{j}(a,b)}
(cd)=
(t_{j+1}(a,b)c)d+
c(t_{j+2}(a,b)d)\eqno(1.5)$$
for any $a,b,c,d\in A$
and for any
$j=0,1,2,$
where the indices are defined
modulo
$3,$
i.e.
$$
t_{j\pm 3}(a,b)=
t_{j}(a,b)\eqno(1.6)$$
for
$j=0,1,2.$
Here,
${\bar Q}\in \ {\rm End}\ A$
for any
$Q\in \ {\rm End}\ A$
is defined by
$$
\overline{ Qa}=
{\bar Q}{\bar a}.\eqno(1.7)$$
We call the unital involutive algebra $A$
satisfying
Eq.(1.5)
to be a pre-structurable algebra.
Moreover,
let
$$
Q(a,b,c):=
t_{0}(a,{\bar b} {\bar c})+
t_{1}(b,{\bar c} {\bar a})+
t_{2}(c,{\bar a} {\bar b})\eqno(1.8)$$
to satisfy 
$Eq (X)$
of [A-F],
which is rewritten as (see [O.1])
$$
Q(a,b,c)=0,\eqno(X)$$
then we call the
pre-structurable algebra
$A$ to be structurable.
\par
We note that these concept are a generalization of 
well known
$" the \ principle \ of \  triality"$,
because the octonion algebra is a structurable algebra.
\par
The main purpose of this
note is
to study properties of
pre-structurable algebra as well as
those satisfied by
$Q(a,b,c)$, and its relation to Lie algebras.
First,
we note the following
Theorem
(see Theorem 3.1 and Lemma 3.6 of [A-F])
\par
{\bf Theorem 1.1}\par
\it
Let $A$ be any unital algebra with involution, and 
 introduce
$A(a,b,c)$
and $B(a,b,c)\in \ {\rm End A}$
by
$$
A(a,b,c)d:=
((da){\bar b})c-
d(a({\bar b}c))\eqno(1.9a)$$
$$
B(a,b,c)d:=
((da){\bar b})c-
d((a{\bar b})c).\eqno(1.9b)$$
Then,
a necessary and sufficient condition that
the unital involutive algebra
$(A,-)$ is
pre-structurable is to have
$$
A(a,b,c)-
A(b,a,c)=
A(c,a,b)-A(c,b,a).\eqno(A)$$
Moreover,
we have
\par
(i)\quad
Eq.(A) implies the validity of
$$
B(a,b,c)-
B(b,a,c)=
B(c,a,b)-
B(c,b,a)\eqno(B)$$
$$
{\rm and}$$
$$
[a,{\bar b},c]-
[b,{\bar a},c]=
[c,{\bar a},b]-
[c,{\bar b}, a]\eqno(A1)$$
\par
(ii)\quad
Eq.(B) implies
$$
[a-{\bar a},b,c]+
[b,a-{\bar a},c]=0\eqno(sk)$$
\par
(iii)\quad
Eq.(sk) implies
$$
[c,{\bar a},b]-
[c,{\bar b},a]=
[c,a,{\bar b}]-
[c,b,{\bar a}]
\eqno(sk\ 1)$$
\par
\quad
where
$[a,b,c]$
is the associator of $A$ defined by
$$
[a,b,c]:=
(ab)c-a(bc).\eqno(1.10)$$
\par
\rm
{\bf Remark 1.2}
\par
Taking the involution of Eq.(sk),
it yields
$$
[a-{\bar a},b,c]=
-[b,a-{\bar a},c]=
[b,c,a-{\bar a}].\eqno(sk)^{'}$$
Also,
Eqs.(A1)
and (sk1)
can be combined to give
$$
[a,{\bar b},c]-
[b,{\bar a},c]=
[c,{\bar a},b]-
[c,{\bar b},a]=
[c,a,{\bar b}]-
[c,b,{\bar a}].\eqno(A.1)^{'}$$
Hereafter in what follows in this note,
$A$ is always designated to
be a pre-structurable algebra over a
field $F$,
unless it is stated otherwise.
\par
{\bf Proposition 1.3}
\par
\it
Under the assumpsion as in above Theorem 1.1,
we have
\par
(i)
$$
\overline{ t_{j}(a,b)}=
t_{3-j}({\bar a},{\bar b})
\eqno(1.11a)$$
\par
(ii)
$$
[t_{j}(a,b),t_{k}(c,d)]=
t_{k}(t_{j-k}(a,b)c,d)+
t_{k}(c,t_{j-k}(a,b)d),
\eqno(1.11b)$$
for any $a,b,c,d\in A$ and for any
$j,k=0,1,2.$
\par
\rm
\vskip 3mm
{\bf Proof}
\par
(i)\quad It is easy to verify  Eq.(1.11a) for
$j=1,$
and $2$,
so that we need to show it for the case
of
$j=0,$
i.e.
$$
\overline{ t_{0}(a,b)}
=t_{0}({\bar a},{\bar b})\eqno(1.12a)$$
or equivalently
$$
t_{0}(a,b)=
l(b{\bar a}-a{\bar b})+
r(b) r({\bar a})-
r(a) r({\bar b}).\eqno(1.12b)$$
This requires the validity of 
$$
\{l(b{\bar a}-a{\bar b})+
r(b)r({\bar a})-
r(a)r({\bar b})\}c$$
$$
=
\{r({\bar a}b-{\bar b}a)+
l(b)l({\bar a})-
l(a)l({\bar b})\}c$$
which is equivalent to the validity of
Eq.(A1),
i.e.
$$
[b,{\bar a},c]-
[a,{\bar b},c]+
[c,{\bar a},b]-
[c,{\bar b},a]=0.$$
\par
(ii)\quad
Eq.(1.11b) follows from Proposition 2.1 of
[O.1]
by considering the conjugate algebra
$A^{*}$ of
$A$ defined by a new bi-linear product
$a*b= \overline{ ab}={\bar b}{\bar a},$
since the condition (C) of [O.1]
 is satisfied in
view of $A$ being unital. $\square$
\par
{\bf Corollary.1.4}
\par
\it        
Under the assumpsion as in above,
let us introduce a triple product in $A$ by
$$
abc:=
t_{0}(a,b)c=
c({\bar a}b-{\bar b}a)+
b({\bar a}c)-
a({\bar b}c).$$
Then,
it defines a Lie triple system,
i.e.
it satisfies
\par
(i)
$$
abc=-bac\eqno(1.13a)$$
(ii)
$$
abc+bca+cab=0\eqno(1.13b)$$
(iii)
$$
ab(cdf)=
(abc)df+
c(abd)f+
cd(abf)\eqno(1.13c)$$
for any $a,b,c,d,f \in A.$
\par
\rm
{\bf Proof}
\par
Eq.(1.13a) is trivial,
while Eq.(1.13b)
follows immediately  from Eq.(1.4).
Finally,
Eq (1.13c) is equivalent to the validity of
$$
[t_{0}(a,b),
t_{0}(c,d)]=
t_{0}(t_{0}(a,b)c,d)+
t_{0}(c,t_{0}(a,b)d)\eqno(1.13d)$$
which is a special case of
Eq.(1.11b)
for
$j=k=0.$  $\square$
\par
{\bf Remark 1.5}
\par
If we set
$$
L(a,b):=
t_{0}(a,b)+
t_{2}({\bar a},{\bar b})$$
then it also satisfies
$$
[L(a,b),L(c,d)]=
L(L(a,b)c,d)+
L(c,L(a,b)d)$$
although we will not go to its detail.
\par
We next set
$$
D(a,b):=
t_{0}(a,b)+
t_{1}(a,b)+
t_{2}(a,b).\eqno(1.14)$$
We then find 
\par
\vskip 3mm
{\bf Proposition 1.6.}
\it
\par
\it
Under the asumpsion as in above,
we have
\par
(i)
$$
\overline{ D(a,b)}=
D({\bar a},{\bar b})=
D(a,b)\eqno(1.15a)$$
(ii)
$$
D(a,b)(cd)=
(D(a,b)c)d+
c(D(a,b)d)\eqno(1.15b)$$
i.e.
$D(a,b)$
is a derivation of $A.$
\par
\noindent
(iii)
$$
[D(a,b),t_{k}(c,d)]=
t_{k}(D(a,b)c,d)+
t_{k}(c,D(a,b)d).\eqno(1.15c)$$
\par
\rm
{\bf Proof}
\par
\noindent
(i)\quad
$\overline{ D(a,b)}=D({\bar a},{\bar b})$
follows immediately from Eq.(1.11a).
\par
We next express
$$
\begin{array}{ll}
D(a,b)c&=
c({\bar a}b-{\bar b}a)+
b({\bar a}c)-
a({\bar b}c)\\
&
+{\bar b}(ac)-{\bar a}(bc)+
(ca){\bar b}-
(cb){\bar a}
\end{array}
$$
by Eqs.(1.2),
so that we calculate
$$
\begin{array}{l}
\{D(a,b)-D({\bar a},{\bar b})\}c\\
=[c,a,{\bar b}]-
[c,b,{\bar a}]+
[c,{\bar b}, a]-
[c,{\bar a},b]
\end{array}
$$
which is identically zero by
Eq.(sk1) to yield
$D(a,b)=D({\bar a},{\bar b}).$
\par
\noindent
(ii)\quad
Summing over $l=0,1,2$
in Eq.(1.5),
it gives
$$
\overline{ D(a,b)}
(cd)=
(D(a,b)c)d+
c(D(a,b)d)$$
and have
Eq.(1.15b)
because of Eq.(1.15a).
\par
\noindent
(iii)\quad
Eq.(1.15c)
follows from Eq.(1.11b)
by summing over $j.$ $\square$
\par
Next,
we assume that underlying field
$F$ to be of
charachteristic
$\not= 2,$
and set
$$
S=\{a|{\bar a}=a,\ 
a\in A\}
\eqno(1.16a)$$
$$
H=\{a|{\bar a}=-a,\ 
a\in A\}\eqno(1.16b)$$
Then $H$ is an algebra with respect to
the anti-commutative 
product
$[a,b]=ab-ba.$
Since Eq.(sk)
implies then
$H$ to be a
generalized alternative nucleus of
$A$,
$H$ is a Malcev algebra with
respect to the product
$[a,b].$
(see [P-S])
\par

{\bf Proposition 1.7}
\par
\it
Let $A$ be a unital involutive algebra over the
field $F$
of charachteristic $\not= 2$,
satisfying Eq.(sk).
If $Dim\ S=1$, then A  is  alternative 
and hence structurable.Moreover it is quadratic, satisfying 
$$a{\bar a}=<a|a>e,\ \ {\bar a}=2<a|e>e-a,$$
for a symmetric bi-linear form $<\cdot |\cdot >$ .
Especially, A is a composition algebra satisfying
 
$$
<ab|ab>=<a|a><b|b>,
$$
$$
<{\bar a}|bc>=<{\bar b}|ca>=<{\bar c}|ab>
$$
.for $ a,b,c\ \in A.$
However $Dim \ A$ needs not be limited to
the canonical value ([S]) of 1,2,4, or 8,
here since $<\cdot |\cdot>$ may be degenerate.  
\par
\rm
{\bf Proof}
\par
Since the unit element $e$ satisfies
$$
[e,a,b]=[a,e,b]=[a,b,e]=0
$$
for any $a,b \ \in A$,
$A$ is clearly alternative if Eq. (sk) holds
valid.
Moreover noting $a{\bar a} \ \in S$ and $a+{\bar a} \ \in S$,
we can write 
$$
a{\bar a}=<a|a>e,\ \ a+{\bar a}=t(a)e
$$
for a linear form $t(a)$ and a symmetric bi-linear form $<\cdot |\cdot>$,
in view of $Dim \ S=1$ Linearizing the
first relation,
it gives
$$
ab+ba=t(a)b+t(b)a-2<a|b>e.
$$
Setting $b=e$,
this yields $t(a)=2<e|a>.$
\par
\noindent
Since $A$ is
a quadratic alteternative algebra, we 
have
$$
(ax){\bar x}={\bar x}(xa)=<x|x>a
$$
for any $a,x \ \in A$,
since 
$$
(ax){\bar x}=[a,x,{\bar x}]+a(x{\bar x})=a(x{\bar x})=<x|x>a$$
and
$$
{\bar x}(xa)= ({\bar x}x)a-[{\bar x},x,a]=<x|x>a
$$
by the alternative law together with quadratic one.
\par
Linearizing them, we have also
$$
(ax){\bar y}+(ay){\bar x}={\bar x}(ya) +{\bar y}(xa)=2<x|y>a.
$$
We can calculate $(ab)(ca)$ in the following two ways.
\par
First, set ${\bar x}=ab$ to compute
$$
(ab)(ca)={\bar x}(ca)=2<x|c>a-{\bar c}(xa)
$$
$$
=
2<\overline{ab}|c>a - {\bar c}\{({\bar b}{\bar a})a\}
$$
$$
=2<ab|{\bar c}>a -{\bar c}\{<a|a>{\bar b}\}
$$
so that 
$$(ab)(ca)=2<{\bar c}|ab>a -<a|a>{\bar c}{\bar b}.\eqno(1.17)
$$
Next, setting ${\bar y}=ca$, 
we similarly proceed 
$$
(ab)(ca)=(ab){\bar y}=2<b|y>a -(ay){\bar b}
$$
$$
=2<b|\overline{ca}>a - \{a({\bar a}\ {\bar c})\}{\bar b}
$$
$$
=
2<{\bar b}|ca>a - \{<a|a>{\bar c}\}{\bar b}
$$
$$
=2<{\bar b}| ca>a - <a|a>{\bar c}{\bar b}.
$$  
Comparing this with Eq. (1.17), we obtain
$$
<{\bar c}|ab>=
<{\bar b}|ca>.
$$
Finally,
$$<ab|ab>=<\overline{{\bar b}\ {\bar a}}|ab>=<{\bar a}|b({\bar b}\ {\bar a})>
$$
$$
=<{\bar a}|<b|b>{\bar a}>=<b|b><{\bar a}|{\bar a}>=<b|b><a|a>.
$$
  This completes the proof of Proposition 1.7. $\square$
 \vskip 2mm 

\par
{\bf Remark 1.8}\par
 
Let us consider the case of $Dim\ A=3$ 
with $S=Fe$ and $H=<f,g>_{span}$.
Then a general solution satisfying Eq. (sk)
is obtained as
$$
f^{2}=\alpha ^{2}e,\ \ g^{2}=\beta ^{2}e,
$$
$$
fg=-\alpha \beta e +\beta f +\alpha g,\ gf=-\alpha \beta e -\beta f - \alpha g \eqno(1.18)
$$
for $\alpha, \ \beta \ \in F,$ satisfying
$$
\alpha ^{2} =-<f|f>,\ \beta ^{2}=-<g|g>,\ \ \alpha \beta =<f|g>=<g|f>,
$$
together with 
$<e|f>=<e|g>=0$
and
$<e|e>=1$. These give 
a composition algebra satisfying  $<ab|ab>=<a|a><b|b>$ with $Dim\ A=3$,
e.g. we have
$$
<f^{2}|f^{2}>=<f|f><f|f>,\ and\ \ <fg|fg>=<f|f><g|g>\ \ etc.
$$
\par
We note that $<\cdot |\cdot >$ is degenerate, since we have  
$$
<(\beta f+\alpha g)|x>=0
$$
for any $x \ \in A.$

This algebra also satisfies 
a linear composition law below.
\par
Let $\phi : A \rightarrow F $ be a linear form defined by
$$
\phi (e)=1,\ \phi (f)=\alpha,\ \phi (g)=-\beta.
$$
Then it satisfies
the linear composition law of
$$
\phi (xy)=\phi (x)\phi (y)
$$
for any $x, y\ \in A.$
\par
We have also found another composition algebra of $Dim\ A=3$
now satisfying $Dim\ S=2$ and $Dim\ H=1$. Let $S=<e,f>_{span}$
and $H=Fg$ for some $f,g \ \in A$.
Suppose that we have
$$
f^{2}={\alpha}^{2} e,\ fg=gf=\alpha g, g^{2}=\lambda (\alpha e+f) \eqno(1.19)
$$
for some $\alpha,\lambda \in F$
Then $A$ is comutative and associative
so that
it is alternative.
Moreover,
$A$ satisfies a composition law $<ab|ab> =
<a|a><b|b>$,
if we introduce a bi-linear symmetric form $<\cdot |\cdot>$ by
$$
<e|e>=1,\ <e|f>=\alpha, \ <e|g>=0,
$$
$$
<f|f>=\alpha^{2}, \ <g|g>=-2\lambda \alpha,\ <f|g>=0.
$$
Indeed, we have
$$
<f^{2}|f^{2}>=<f|f><f|f>\ \ and\ <fg|fg>=<f|f><g|g>, \ etc.
$$
Note that $
<\cdot |\cdot>$ is degenerate, satisfying
$$
<(f-\alpha e)|x>=0
$$
for any $x\ \in A$.
\par
Moreover, defining a linear form $\phi :A \rightarrow F$
by $$
\phi (e)=1,\ \phi (f)=\alpha,\ (\phi (g))^{2}=2\lambda \alpha,
$$
it also satisfies a linear composition law
$$
\phi (xy)=\phi (x)\phi (y).
$$

\par
{\bf Examples 1.9}
Let
$$
A=\left(
\begin{array}{ll}
F&B\\
B&F
\end{array}
\right)\eqno(1.20)$$
be the Zorn's vector matrix algebra,
where
$B$ is a involutive algebras over the field
$F$ with a bi-linear form
$(\cdot|\cdot)$.
We introduce a bi-linear product
in $A$ by
$$
X_{1}X_{2}:=
\left(
\begin{array}{ll}
\alpha_{1}&x_{1}\\
y_{1}&\beta_{1}
\end{array}
\right)
\left(
\begin{array}{ll}
\alpha_{2}&x_{2}\\
y_{2}&\beta_{2}
\end{array}
\right):=
\left(
\begin{array}{ll}
\alpha_{1}\alpha_{2}+(x_{1}|y_{2}),&
\alpha_{1}x_{2}+\beta_{2}x_{1}+
y_{1}y_{2}\\
\alpha_{2}y_{1}+\beta_{1}y_{2}+
x_{1}x_{2},&
\beta_{1}\beta_{2}+
(y_{1}|x_{2})
\end{array}\right)
\eqno(1.21)$$
for 
$\alpha_{\bar j},\ \beta_{\bar j}\in F$
and
$x_{\bar j},\ y_{\bar j}\in B$
for $j=1,2.$
Then $A$ is unital with 
the unit element
$$
E=
\left(
\begin{array}{ll}
1&0\\
0&1
\end{array}\right)\eqno(1.22)
$$
and is involutive with involution map
$$
\left(
\begin{array}{ll}
\alpha&x\\
y&\beta
\end{array}
\right)
\rightarrow
\overline{
\left(
\begin{array}{ll}
\alpha&x\\
y&\beta
\end{array}
\right)}=
\left(
\begin{array}{ll}
\beta&{\bar x}\\
{\bar y}&\alpha
\end{array}
\right),
\eqno(1.23)
$$
provided that
$(\cdot|\cdot)$ satisfies
$$
({\bar y}|{\bar x})=
(x|y).\eqno(1.24)$$
Some case of $B$ being commutative with
${\bar x}=x$
has been studied
in ([A-F],[O.1]
and [Kam.1]),
where $B$ is an admissible cubic
 algebra
([E-O.1]).
Here,
we will consider the another case of $B$ 
being anti-commutative i.e.
$xy=-yx$ with
${\bar x}=-x$.
In that case,
we see
$Dim\ S=1,$
since only even 
element of $A$ is the unit element $E$.
The condition that $A$ satisfies
Eq.(sk) is obtained
if and only if we have
$$
x(yz)=
(y|x)z-
(x|z)y\eqno(1.25a)$$
$$
(x|yz)=
(y|zx)=
(z|xy).\eqno(1.25b)$$
Then, by Proposition 1.7,
$A$ is an alternative algebra  and its anticommutative
algebra given by
$[X,Y]=
XY-YX$
is  a Malcev algebra.
In other words,
$A$ is a flexible
Malcev admissible 
algebra.(see [M]), and
 $A$ is quadratic.
Now    
consider the case of $Dim\ B =3,$ where
$B$ is generated by three
element 
$e_{1},e_{2}$
and $e_{3}$
satisfying
$$
\begin{array}{l}
e_{1}e_{2}=-e_{2}e_{1}=e_{3},\\
e_{2}e_{3}=-e_{3}e_{2}=e_{1},\\
e_{3}e_{1}=-e_{1}e_{3}=e_{2},\\
e_{1}e_{1}=e_{2}e_{2}=e_{3}e_{3}=0,
\end{array}
$$
with
$(e_{i}|e_{j})=
\delta_{ij}
\ 
(i,j,=
1,2,3).$
Then,
these satisfy Eqs.(1.25)
and $A$ is an octonion algebra since this case
corresponds to the original Zorn
vector matrix algebra
[S].
Similarly,
the case of 
$Dim\ B=1$ with
$B=Fe_{0}$
satisfying
$e_{0}e_{0}=0$
and
$(e_{0}|e_{0})=1$
gives the quaternion algebra.
\par
However the case of
$Dim\ B=2$
with 
$B=<e_{1},e_{2}>_{span}$
for [A,A]
will yield a
$5$-dimensional Malcev algebra of
Kuzmin type
([Ku],[M]),
when we impose
$$
\begin{array}{ll}
e_{1}e_{2}=-e_{2}e_{1}=e_{1},&
e_{1}e_{1}=e_{2}e_{2}=0,\\
(e_{1}|e_{1})=
(e_{1}|e_{2})=0,&
(e_{2}|e_{2})=1.
\end{array}
$$
Note then that
[A,A]
is a non-Lie Malcev algebra.
\par

We introduce a bi-linear symmetric form in $A$ by
$$
<X_{1}|X_{2}>=1/2\{ \alpha_{1}\beta_{2}+\beta_{1}\alpha_{2}-(x_{1}|y_{2})-(y_{1}|x_{2})\},
$$
it satisfies the quadratic law
$$
X_{1}{\bar X_{1}}
=<X_{1}|X_{1}>E
$$
and 
$$
<{\bar X_{1}}|X_{2}X_{3}>
 =<{\bar X_{2}}|X_{3}X_{1}>
 =<{\bar X_{3}}|X_{1}X_{2}>.
 $$Moreover since  $A$
is also alternative,
then $A$ satisfies
also the composition law;
$$
<X_{1}X_{2}|X_{1}X_{2}> =<X_{1}|X_{1}><X_{2}|X_{2}>.
$$ 
However, we note that $<X_{1}|X_{2}>$ is degenerate,whenever
$(x_{1}|x_{2})$ in $B$
is degenerate .
Therefore existence of the composition algebra $A$ with $Dim\ A =6$
corresponding to the case of $Dim \ B=2$ does not contradict the Hurwitz theorem([S]). 
\par

Finally,
if we suppose $B$ to be a trivial algebra with
$xy=0$
and $(x|y)=0$
for all
$x,y\in B,$
then the  algebra $A$ becomes a 
commutative associative algebra.
\par
\vskip 3mm
{\bf II The property of  $Q(a,b,c)$}
\par
\vskip 3mm
Here, we assume $A$ to be a pre-structurable algebra and
we will discuss 
some properties of
$Q(a,b,c).$
\par
\vskip 3mm
{\bf Theorem 2.1}
\par
\it
Let $A$ be a pre-structurable algebra.
Then,
we have
\par
(i)\ $Q(a,b,c)d$ is totally symmetric
in $a,b,c,d\in A.$
\par
(ii)\ $Q(a,b,c)d$ is identically zero,
if at least one of
$a,b,c,d$
is the identity element $e.$
\par
(iii)\ Suppose the underlying field $F$
to be of characteristic
$\not= 2,$
then
$Q(a,b,c)d$
is identically 
\par
zero,
if at least one of 
$a,b,c$
and $d$ is a element of
$H.$
\par
(iv)\ $\overline{ Q(a,b,c)}=
Q({\bar a},{\bar b},{\bar c})
=Q(a,b,c)$
is a derivation of $A.$
\par
(v)\ 
$3Q(a,b,c)=D(a,{\bar b}{\bar c})+
D(b,{\bar c}\ {\bar a})+
D(c,{\bar a}{\bar b}).$
\par
\rm
For a proof of this theorem,
we start with the following Lemmas:
\par
\vskip 3mm
{\bf Lemma 2.2}
\par
\it
Under the asumpsion as in above,
$$
Q(a,b,c)d=Q(d,b,c)a\eqno(2.1)$$
is symmetric in $a$ and
$d.$
\par
\vskip 3mm
{\bf Proof}
\par
\rm
We proceed as follows:
Since
$$
Q(a,b,c)d=
(t_{0}(a,{\bar b}{\bar c})+
t_{1}(b,{\bar c}\ {\bar a})+
t_{2}(c,{\bar a}{\bar b}))d,$$
we note
$$
\begin{array}{l}
Q(a,b,c)d-Q(d,b,c)a\\
=\{t_{0}(a,{\bar b}{\bar c})d-
t_{0}(d,{\bar b}{\bar c})a\}+
\{t_{1}(b,{\bar c}\  {\bar a})+
t_{2}(c,{\bar a}{\bar b})\}d\\
-\{t_{1}(b,{\bar c}{\bar d})+
t_{2}(c,{\bar d}\   {\bar b})\}a.
\end{array}$$
However, we calculate
$$
\begin{array}{l}
t_{0}(a,{\bar b}{\bar c})d-
t_{0}(d,{\bar b}{\bar c})a=
-t_{0}({\bar b}{\bar c},a)d-
t_{0}(d,{\bar b}{\bar c})a\\
=
t_{0}(a,d)({\bar b}{\bar c})=
\overline{ t_{0}({\bar a},{\bar d})}
({\bar b}{\bar c})=
(t_{1}({\bar a},{\bar d}){\bar b}){\bar c}+
{\bar b}
(t_{2}({\bar a},{\bar d}){\bar c})
\end{array}
$$
so that
$$
\begin{array}{l}
Q(a,b,c)d-Q(d,b,c)a\\
=(t_{1}({\bar a},{\bar d}){\bar b}){\bar c}+
{\bar b}(t_{2}({\bar a},{\bar d}){\bar c})+
t_{1}(b,{\bar c}\ {\bar a})d+t_{2}(c,{\bar a}{\bar b})d\\
-t_{1}(b,{\bar c}{\bar d})a-
t_{2}(c,{\bar d}\ {\bar b})a\\
=\{(l(d)l({\bar a})-
l(a)l({\bar d})){\bar b}\}{\bar c}+
{\bar b}\{r(d)r({\bar a})
-r(a)r({\bar d})\}{\bar c}\\
+\{l(ac)l(b)-
l({\bar b})l({\bar c}{\bar a})\}d+
\{r(ba)r(c)-
r({\bar c})
r({\bar a}{\bar b})\}d\\
-\{l(dc)l(b)-l({\bar b})l({\bar c}{\bar d})\}a-
\{r(bd)r(c)-
r({\bar c})r({\bar d}\  {\bar b})\}a\\
=\{d({\bar a}{\bar b})-
a({\bar d}\ {\bar b})\}{\bar c}+
{\bar b}\{({\bar c}\ {\bar a})d-
({\bar c} {\bar d})a\}\\
+
(ac)(bd)-
{\bar b}\{({\bar c}\  {\bar a})d\}+
(dc)(ba)-
\{d({\bar a} {\bar b})
\}{\bar c}\\
-(dc)(ba)+
{\bar b}\{({\bar c} {\bar d})a\}-
(ac)(bd)+
\{a({\bar d}\  {\bar b})\}{\bar c}\\
=0
\end{array}
$$
identically,
proving Eq.(2.1). $\square$
\par
\vskip 3mm
{\bf Lemma 2.3}
\par
\it
Under the assumpsion as in above,
we have
$$
Q(a,b,c)e=Q(e,b,c)a=0$$
for any $a,b,c\in A,$
where $e$ is the unit element of
$A.$
\par
\vskip 3mm
{\bf Proof}
\rm
We caluculate
$$
\begin{array}{l}
Q(e,b,c)a=\{t_{0}(e,{\bar b}{\bar c})
+t_{1}(b,{\bar c}{\bar e})+
t_{2}(c,{\bar e}\ {\bar b})\}a\\
=
\{r({\bar e}({\bar b}{\bar c})-
(cb)e)+
l({\bar b}{\bar c})l({\bar e})-
l(e)l(cb)\\
+
l(c)l(b)-l({\bar b})l({\bar c})+
r(b)r(c)-r({\bar c})
r({\bar b})\}a\\
=
a({\bar b}{\bar c}-cb)+
({\bar b}{\bar c})a-
(cb)a\\
+
c(ba)-{\bar b}
({\bar c}a)+
(ac)b-
(a{\bar b}){\bar c}\\
=-
[a,{\bar b},{\bar c}]+
[a,c,b]+
[{\bar b},{\bar c},a]-
[c,b,a]=0
\end{array}
$$
by Eq.$(A1)'$
if we let
$a\rightarrow {\bar b}\rightarrow
{\bar c}
\rightarrow {\bar a}.$
$\square$
\par
\vskip 3mm
{\bf Lemma 2.4}
\par
\it
Under the assumpsion as in above, we have
\par
(i)
$$
\overline{ Q(a,b,c)}=Q({\bar a},{\bar c},{\bar b})\eqno(2.2)$$
\par
(ii)
$$
\overline{ Q(a,b,c)}
(df)=
\{Q(c,a,b)d\}f+
d\{Q(b,c,a)f\}.\eqno(2.3)
$$
\vskip 3mm
{\bf Proof}
\par
\rm
(i)
$$
\begin{array}{l}
\overline{ Q(a,b,c)}=
\overline{ t_{0}(a,{\bar b} {\bar c})}
+
\overline{ t_{1}(b,{\bar c} {\bar a})}
+
\overline {t_{2}(c,{\bar a} {\bar b})}\\
=
t_{0}({\bar a},cb)+
t_{2}({\bar b},ac)+
t_{1}({\bar c},ba)=
Q({\bar a},{\bar c},{\bar b})
\end{array}
$$
by Eq.(1.11a).
\par
(ii)
\par
By Eq.(1.5),
we have
$$
\begin{array}{l}
\overline{ t_{0}(a,{\bar b}{\bar c})}
(df)=
\{t_{1}(a,{\bar b}{\bar c})d\}f+
d\{t_{2}(a,{\bar b}{\bar c})f\}\\
\overline{ t_{1}(b,{\bar c}\  {\bar a})}(df)=
\{t_{2}(b,{\bar c}\ {\bar a})d\}f+
d\{t_{0}(b,{\bar c}\  {\bar a})f\}\\
\overline{ t_{2}(c,{\bar a} {\bar b})}
(df)=
\{t_{0}(c,{\bar a}{\bar b})d\}f+
d\{t_{1}(c,{\bar a}{\bar b})f\}.
\end{array}
$$
Summing over all these relations,
this give Eq.(2.3).$\square$
\par
We are now in position to  prove Theorem 2.1.
We first set
$d=e$
or $f=e$ and then letting $f  \rightarrow d $  in Eq.(2.3).
\par
We find
$$
\overline{ Q(a,b,c)}d=
Q(c,a,b)d=
Q(b,c,a)d
\eqno(2.4)$$
which shows the
$Q(c,a,b)$ is cyclic invariant
under
$ a\rightarrow b\rightarrow c\rightarrow a,$
and hence,
it also gives
$$
\overline{ Q(a,b,c)}d=
Q(a,b,c)d.,$$
i.e.
$$
\overline{ Q(a,b,c)}=
Q(a,b,c).\eqno(2.5)$$
Further
$Q(a,b,c)d=
Q(d,b,c)a$
by Lemma 2.2,
so that
$$
\begin{array}{l}
Q(a,b,c)d=
Q(d,b,c)a=
Q(b,c,d)a=
Q(a,c,d)b\\
=Q(d,a,c)b=Q(b,a,c)d.
\end{array}
$$
Therefore,
we have also
$Q(a,b,c)=Q(b,a,c)$
so that
$Q(a,b,c)$
is totally symmetric in
$a,b,c\in A.$
Especially, this implies
$Q(a,b,c)d$ to be totally symmetric in
$a,b,c$ and $d.,$
proving (i) of Theorem 2.1.
Then together with Lemma 2.2,
it also proves (ii) of
Theorem 2.1.
In order to show
(iii),
we first note the validity of
$$
\overline{ Q(a,b,c)}=
Q({\bar a},{\bar b},{\bar c})=
Q(a,b,c)
\eqno(2.6)
$$
because of Eqs.(2.2) and
(2.5).
This gives (iv) of the Theorem.
\par
In order to avoid possible confusion,
we label elements of $S$ and $H,$
respectively
as $a_{0}$
and $a_{1},$
or $b_{0}$ and
$b_{1}$ etc.
so that
$$
{\bar a_{0}}=a_{0},\ {\bar a_{1}}=-a_{1}\ {\rm etc.}
$$
Then,
Eq.(2.6) implies immediately
$$Q(a_{0},b_{0},c_{1})=0=
Q(a_{1},b_{1},c_{1}).$$
\par
Note that we have
$$
Q(a_{0},b_{1},c_{1}) d_{1}=
Q(d_{1},b_{1},c_{1}) a_{0}=0$$
i.e.,
$$
Q(a_{0},b_{1},c_{1})d_{1}=0.$$
However, we have also
$$
Q(a_{0},b_{1},c_{1})d_{0}=
Q(a_{0},d_{0},c_{1})b_{1}=0$$
so that
$Q(a_{0},b_{1},c_{1})=0.$
This proves the statement
(iii) of Theorem 2.1.
Finally,
$$
Q(a,b,c)=
t_{0}(c,{\bar b}{\bar c})+
t_{1}(b,{\bar c}{\bar a})+
t_{2}(c,{\bar a}{\bar b}).$$
Letting $a\rightarrow b\rightarrow c\rightarrow a,$
and adding all of the resulting relation
we obtain
$$
Q(a,b,c)+Q(b,c,a)+
Q(c,a,b)=
D(a,{\bar b}{\bar c})+
D(b,{\bar c}{\bar a})+
D(c,{\bar a}{\bar b})$$
so that
$$
3Q(a,b,c)=
D(a,{\bar b}{\bar c})+
D(b,{\bar c}{\bar a})+
D(c,{\bar a}{\bar b}).\eqno(2.7)$$
This completes the proof of
Theorem 2.1.$\square$
\par
\vskip 3mm

{\bf Proposition 2.5}
\par
\it
Under the assumpsion as in above,
if $a={\bar a}\in S,$
then we have
$$
\begin{array}{ll}
Q(a,a,a)a&
=[a,aa^{2}]+3\{a^{2}a^{2}-a(a^{2}a)\}\\
&
=[a^{2}a,a]+
3\{a^{2}a^{2}-(aa^{2})a\}.
\end{array}
$$
\par
\vskip 3mm
{\bf Proof}

\rm
\par
We calculate
$$
\begin{array}{ll}
Q(a,a,a)&
=t_{0}(a,{\bar a}{\bar a})+
t_{1}(a,{\bar a}{\bar a})+
t_{2}(a,{\bar a}{\bar a})\\
&
=t_{0}(a,a^{2})+t_{1}(a,a^{2})+
t_{2}(a,a^{2})\\
&
=r(aa^{2}-a^{2}a)+
l(a^{2})l(a)-
l(a)l(a^{2})+
l(a^{2})l(a)-l(a)l(a^{2})\\
&
+r(a^{2})r(a)-
r(a)r(a^{2})
\end{array}
$$
so that
$$
\begin{array}{ll}
Q(a,a,a)a&
=a(aa^{2}-a^{2}a)+
2\{a^{2}a^{2}-
a(a^{2}a)\}+
a^{2}a^{2}-
(aa^{2})a\\
&=
[a,aa^{2}]+
3\{a^{2}a^{2}-a(a^{2}a)\}.
\end{array}
$$
\par
Taking the involution of this relation,
we have also
$$
Q(a,a,a)a=
[a^{2}a,a]+3\{a^{2}a^{2}-
(aa^{2})a\},
$$
these prove the Proposition.$\square$
\par
\vskip 3mm
{\bf Remark 2.6 (a)}
\par
If $A$ is a pre-structurable algebra 
over the field $F$ of charachtristic
$\not= 2,$
and $\not= 3,$
then Theorem 2.1 and Proposition 2.5
imply that
$A$ is structurable,
provided that we have
$aa^{2}=a^{2}a(:=a^{3}$
and $a^{2}a^{2}=aa^{3}(=a^{3}a)$
for $a\in S,$
since then
$Q(a,b,c)=0$
for any $a,b,c\in A.$
\par
Thus if $A$ is power-associative and pre-structurable algebra, then $A$ is structurable.
\par
{\bf Remark 2.6 (b)}\quad
\par
By means of relations (2.3) and(2.6),
we note that
$$
Q(a,b,c)\ is\  a\  \ derivation \ of \ pre-structurable \ algebra.
 \  A.$$
\vskip 5mm
{\bf Proposition 2.7}
\par
\it 
Let $A$ be a pre-structurable algebra and set
$A_{0}=\{x|x\in A,$
and $Q(a,b,c)x=0$
for any $a,b,c\in A\}.$
Then,$A_{0}$ is a structurable algebra.
\par
\vskip 3mm
{\bf Proof}
\par
\rm
First,
let us show that
$A_{0}$ is a sub-algebra of $A$,
since we calculate
$$
Q(a,b,c)(xy)=
(Q(a,b,c)x)y+
x(Q(a,b,c)y)=0$$
for any
$x,y\in A_{0}$ so that
$xy\in A_{0}$.
Moreover,
$e\in A_{0}$ by
Theorem 2.1.
Also,
if $x\in A_{0},$ then
${\bar x}\in A_{0}$
since
$$
0=
\overline{ Q(a,b,c)x}=
\overline{ Q(a,b,c)}{\bar x}=
Q(a,b,c){\bar x}$$
by (iv) of Theorem 2.1.
Then,
$A_{0}$ is involutive.
Further,
we note
$d_{j}(x,y)\in {\rm End}\ A_{0}$
for
$x,y\in A_{0},$
proving that
$A_{0}$ is a pre-structurable algebra.
Finally
since $Q(x,y,z)=0$
restricted to $A_{0}$ and,
this proves $A_{0}$ to be structurable. $\square$
\par
By Theorem 2.1,
$A_{0}$ always
contains a structurable subalgebra generated by
$e$ and elements of $H$.
It is plausible that we may have
$A_{0}=A$
if $A$ is a simple algebra.
However, we could neither prove nor disprove 
such a conjecture.
\par
We will now prove the converse statement of Theorem 2.1.
\par
{\bf Theorem 2.8}\it
\par 
Let A be a unital involutive algebra satisfying 
\par
(i)  $Q(a,b,c)d$ is totally symmetric in $a,b,c,d$,
\par
(ii) $\overline{Q(a,b,c)} = Q(a,b,c)$,
\par
(iii) $Q(a,b,c)=0 $ whenever at least one of $a,b,c \in A$ is 
a element of $H$, 
\par
(iv) the validity of Eq. (sk).
\par
Then A is a pre-structurable algebra.
\par
\noindent
Alternatively any unital involutive algebra A is pre-structurable, if Eq. (sk) holds valid and 
we have
$$
Q(a,b,c)=B(b,a,c)-C(a,b,c)-C(c,b,a)-C^{' }(c,a,b) \eqno(2.8)
$$
being totally symmetric in $a,b,c\in A.$
\par
\noindent
Here $C(a,b,c)$ and
$ C^{'}(a,b,c) \in End\ A$
are defined by (see [A-F])
$$
C(a,b,c)d=\{a({\bar b}\ {\bar d})\}c- (a{\bar b})({\bar d}c) \eqno(2.9a)
$$
$$
C^{'}(a,b,c)d=C(a,b,c){\bar d}=\{a({\bar b}\ d)\}c-(a{\bar b})(dc) \eqno(2.9b)
$$ 
\rm
{\bf Proof}
\par
We first show that the validity of $(i)-(iii)$ will give Eq. (2.8).
\par
For this, we note 
$$
Q({\bar a}, b,c)=Q(a,{\bar b},c)=Q(a,b,{\bar c})=Q(a,b,c) \eqno(2.10)
$$

since $a-{\bar a}\in H$ etc. Moreover,
we have also 
$$
\overline{Q(a,b,c)d}=Q(a,b,c)d=Q(a,b,c){\bar d} \eqno(2.11)$$
since we calulate
$$
Q(a,b,c)(d-{\bar d})=Q(a,b,d-{\bar d})c=0.
$$

Next Eqs. (1.2) and (1.8) give 
$$
Q(a,b,c){\bar d}=
t_{0}(a, {\bar b}\ {\bar c}){\bar d}
+t_{1}(b,  {\bar c}\ {\bar a}){\bar d}
+t_{2}(c, {\bar a}\ {\bar b}){\bar d}$$
$$
={\bar d}\{{\bar a}({\bar b} {\bar c})-
(cb)a\}+
({\bar b}{\bar c})({\bar a}{\bar d})-a\{(cb){\bar d}\}$$
$$
+(ac)(b{\bar d})-{\bar b}\{({\bar c}\ {\bar a}){\bar d}\}+
({\bar d}c)(ba)-
\{{\bar d}({\bar a}{\bar b})\}{\bar c}.$$
Taking the involution of this equation
and noting Eq.(2.11),
we find
$$
Q(a,b,c)d=\overline{ Q(a,b,c){\bar d}}$$
$$
=\{(cb)a\}d
-\{ {\bar a} ({\bar b}{\bar c}) \} d
+(da)(cb)
-\{ d({\bar b}{\bar c}) \} {\bar a} $$
$$
+(d{\bar b})({\bar c}\ {\bar a})
-\{d(ac)\}b
+({\bar a} {\bar b}) ({\bar c}d)
-c \{ (ba)d \} $$ 
$$
=\{(cb)a\}d
-c\{(ba)d\}
+(da)(cb)
-\{d(ac)\}b$$
$$
+({\bar a}{\bar b})({\bar c}d)
-\{{\bar a}({\bar b}{\bar c})\}d
+(d{\bar b})({\bar c}\ {\bar a})
-\{d({\bar b}{\bar c})\}{\bar a}$$
$$
=B(b,{\bar a},d)c
-C({\bar a},b,d)c
-C(d,b,{\bar a})c
-C^{'}( d,{\bar a},b)c
$$
However, since
$Q(a,b,c)d=Q(a,b,d)c,$
this is rewritten as
$$
Q(a,b,d)
=B(b,{\bar a},d)-C({\bar a},b,d)
-C(d,b,{\bar a})-C^{'}(d,{\bar a},b)$$
which is Eq.(2.8)
if we let 
$a\rightarrow {\bar a}$
and 
$c\leftrightarrow d$
and noting
$Q({\bar a},b,c)=Q(a,b,c).$
\par
We next prove that
the validity of
Eq.(sk) and Eq.(2.8)
will lead to $A$ to be pre-structurable 
as follows:
\par
Since
$Q(a,b,c)=Q(c,b,a),$
Eq.(2.8)
yield
$$
B(b,a,c)
-B(b,c,a)
=C^{'}(c,a,b)
-C^{'}(a,c,b).\eqno(2.12)$$
Similarly,
we calcultate
$$
\begin{array}{ll}
0&=Q(a,b,c)-Q(b,a,c)\\
&=B(b,a,c)-C(a,b,c)-C(c,b,a)-C^{'}(c,a,b)\\
&-B(a,b,c)+
C(b,a,c)+
C(c,a,b)+
C^{'}(c,b,a).
\end{array}
$$
Then,
together with
Eq.(2.12)
for
$b \leftrightarrow c$,
we obtain
$$
\{B(b,a,c)-B(a,b,c)+B(c,a,b)-B(c,b,a)\}d
$$
$$
=\{C^{'}(b,a,c)-C^{'}(a,b,c)+
C^{'}(c,a,b)-C^{'}(c,b,a)\}d
$$
$$
-\{C(b,a,c)-C(a,b,c)+
C(c,a,b)-C(c.b.a)\}d$$
$$
=-\{C(b,a,c)-C(a,b,c)+C(c,a,b)-
C(c,b,a)\}(d-{\bar d})\eqno(2.13)$$
by recalling
$C^{'}(a,b,c)d=C(a,b,c){\bar d}.$
\par
We now follow the
reasoning
given in section 5 of [A-F]
to obtain
$$
C(a,b,c)(d-{\bar d})=
B(a,b,c)(d-{\bar d})\eqno(2.14)$$
if Eq.(sk) holds,
since we calculate
(with $s=d-{\bar d}$),
$$
\begin{array}{l}
C(a,b,c)s=
\{a({\bar b}{\bar s})\}c-(a{\bar b})({\bar s}c)=
-\{a({\bar b} s)\}c+(a{\bar b})(sc)\\
=\{[a,{\bar b},s]-(a{\bar b})s\}c+
(a{\bar b})(sc)=
[a,{\bar b},s]c-
[a{\bar b},s,c]\\
=[s,a,{\bar b}]c+
[s,a{\bar b},c]=
\{(sa){\bar b}-s(a{\bar b})\}c+
\{s(a{\bar b})\}c-
s\{(a{\bar b})c\}\\
=\{(sa){\bar b}\}c-s\{(a{\bar b})c\}=
B(a,b,c)s.
\end{array}
$$
Then,
Eq.(2.13) becomes
$$
\{B(b,a,c)
-B(a,b,c)+
B(c,a,b)-B(c,b,a)\}d$$
$$
=-\{B(b,a,c)-
B(a,b,c)+
B(c,a,b)-
B(c,b,a)\}
(d-{\bar d}).\eqno(2.15)$$
From Eq.(2.12),
we have also
$$
\begin{array}{l}
\{B(b,a,c)-B(b,c,a)\}
(d-{\bar d})=
\{C^{'}(c,a,b)-
C^{'}(a,c,b)\}
(d-{\bar d})\\
=-\{C(c,a,b)-
C(a,c,b)\}
(d-{\bar d})=
-\{B(c,a,b)-
B(a,c,b)\}(d-{\bar d})
\end{array}
$$
and hence
$$
\{B(b,a,c)-B(b,c,a)+
B(c,a,b)-B(a,c,b)\}(d-{\bar d})=0.$$
Letting $b\leftrightarrow c$,
this together with
Eq.(2.15) gives then
$$
\{B(b,a,c)-B(a,b,c)+
B(c,a,b)-B(c,b,a)\}d=0$$
i.e.,
the validity of Eq.(B).
\par
Finally,
by Eqs.(1.9),
we note
$$
A(a,b,c)=B(a,b,c)+
r([a,{\bar b},c]),$$
where $r(x)d:=dx$
so that
$$
\begin{array}{l}
A(b,a,c)-A(a,b,c)+
A(c,a,b)-A(c,b,a)\\
=B(b,a,c)-B(a,b,c)+
B(c,a,b)-B(c,b,a)\\
+r([b,{\bar a},c]-[a,{\bar b},c]+
[c,{\bar a},b]-[c,{\bar b},a])=0,
\end{array}
$$
provided that
Eq.(A.1) holds in addition.
However,
Eq.(A.1) is also a consequence of
Eq.(2.12) for
$b=e$,
since
$$
\{B(e,a,c)-
B(e,c,a)\}d=
\{C^{'}(c,a,e)-
C^{'}(a,c,e)\}d$$
gives
$$
[d,{\bar a},c]-
[d,{\bar c},a]=
[a,{\bar c},d]-
[c,{\bar a},d]$$
which is Eq.(A.1)
if we let
$d\rightarrow c\rightarrow b.$
Therefore,
Eq.(A) holds valid,
proving $A$ to be
pre-structurable.
This completes the proof of Theorem 2.8. $\square$
\par
\vskip 3mm
The special case of
$Q(a,b,c)=0$
in Theorem 2.8 reproduces
(iii) of Theorem 5.5 of
[A-F]:
\par
\vskip 3mm
{\bf Corollary 2.9}
\par
\it
A necessary and sufficient condition for a unital involutive 
algebra
to be structurable is the validity of Eq.(sk) and Eq.(X).
\par
\rm
We also note the following Proposition
(see Corollary 3.6 of [O.1]):
\par
\vskip 3mm
{\bf Proposition 2.10}
\par
\it
Let $A$ be a unital involutive algebra possessting
a symmetric bi-linear
non-degenerate form
$<\cdot|\cdot>$
satisfying
$$
<{\bar a}|bc>=
<{\bar b}|ca>=<{\bar c}|ab>.$$
Then,
Eq.(A) is equivalent to
Eq.(X).
Especially,
we have
\par
(i)\ 
Any pre-structurable algebra is automatically structurable.
\par
(ii)\ 
If $Q(a,b,c)=0,$
then $A$ is structurable.
\par
\vskip 3mm
\rm
{\bf Remark 2.11}
\par
Many interesting unital involutive algebras containing Jordan and 
alternative algebras are structurable (see [A-F]). It is rather hard to find
example of a pre-structurable but not structurable algebras.
\par
Let $A$ be a
commutative algebra with $Dim\ A=Dim\ S=3$
so that ${\bar a}=a$ for any $a \in A$ and Eq.(sk) is trivially
satisfied. Let $A=<e,f,g>_{span}$ with $e$ being the unit element. Suppose hat we have
$$
ff=fg=gf=0,\ \ gg=\alpha e+\beta f
$$
for $\alpha , \ \beta\ \in F.$
Then we can readily
verify that we have
$Q(a,b,c)d=0$ for any
$a,b,c,d $ assuming values of $e,f,$ and $ g$ except for the case of
$$Q(g,g,g)g= (3\alpha \beta) f.
$$
Therefore, $A$ is pre-structurable but not structurable
by Theorem 2.8, provided that $3\alpha \beta \not=0.$
Note that the case of $\alpha \beta =0$ (or more strongly $\alpha =0$)
corresponds
to
$A$ being a Jordan (or associative) algebra which is structurable.
Note also that the present algebra is not simple since
$B=Ff$ is a ideal of $A$.
\vskip 3mm

\rm
{\bf Remark 2.12}
\par
If $A$ is a structurable algebra,
then
$D(a,b)$
given by
Eq.(1.14)
is a derivation
of $A$ satisfying
$$
D(a,{\bar b} {\bar c})+
D(b,{\bar c} \  {\bar a})+
D(c,{\bar a} {\bar b})=0\eqno(2.16)$$
by Eq.(2.7).
If we set
$$
D_{0}(a,b)=
D(a,{\bar b})\eqno(2.17)$$
then it satisfies
\par
(i)
$$
D_{0}(a,b)=
-D_{0}(b,a)=
D_{0}({\bar a},{\bar b})=
\overline{ D_{0}(a,b)}\eqno(2.18)$$
\par
(ii)
$$
D_{0}(a,b)\ 
{\rm is\ a\ derivation\ of}
A\eqno(2.19a)$$
\par
(iii)
$$
D_{0}(a,bc)+
D_{0}(b,ca)+
D_{0}(c,ab)=0.
\eqno(2.19b)$$
Any algebra $A$ possessing
a non-zero
$D_{0}(a,b)$
satisfying
$D_{0}(a,b)=-D_{0}(b,a)$
as well as Eq.(2.19a)
and
(2.19b)
has been called in [Kam.2]
to be a generalized
structurable algebra.
Therefore,
any structurable algebra is also a 
generalized structurable algebra,
provided that
$D(a,b)\not= 0.$
Note that there exists a
structurable algebra such that 
we have
$D(a,b)=0$
identically as in
Example 25.3 of [O.2].
\par
\vskip 3mm
{\bf III\ 
Lie Algebras with Triality }
\par
\rm
In Corollary 1.4,
we have seen that we can introduce
a Lie triple system for any
pre-structurable algebra
and hence we can construct a Lie algebra in a canonial 
way as follows.
\par
Let
$$
L_{0}=\rho_{0}(A)\oplus T_{0}(A,A)\eqno(3.1)$$
where
$\rho_{0}(A)$
is a copy of $A$
itself
and
$T_{0}(a,b)$,
for $a,b\in A$
is an analogue
(or generalization)
of $t_{0}(a,b)$.
If we wish,
we may identify
$T_{0}(a,b)$ as
$t_{0}(a,b).$
Then,
supposing
commutation 
relations;
\par
(i)
$$
[T_{0}(a,b),T_{0}(c,d)]=
T_{0}(t_{0}(a,b)c,d)+
T_{0}(c,t_{0}(a,b)d)\eqno(3.2a)$$
\par
(ii)
$$
[T_{0}(a,b),\rho_{0}(c)]=
\rho_{0}(t_{0}(a,b)c)\eqno(3.2b)$$
\par
(iii)
$$
[\rho_{0}(a),\rho_{0}(b)]=
T_{0}(a,b)\eqno(3.2c)$$
for
$a,b,c,d\in A,$
$L_{0}$
becomes a Lie algebra as
we may easily verify.
Note that
Eq (3.2a)
is an analogue
of Eq.(1.13d).
\par
A extra advantage of $A$
being structurable is that
we can further enlarge the Lie algebra  
by
utilizing
Eq.(1.11b)
for any
$j,k=0,1,2$
as follows:
\par
Let $\rho_{j}(A)$
for $j=0,1,2$
be 3 copies of $A$.
Moreover,
we introduce three
unspecified
symbols
$T_{j}(a,b)$
for
$j=0,1,2$
and for
$a,b\in A,$
which may be regarded as
a generalization of $t_{j}(a,b).$
If we wish,
we may identify
$T_{j}(a,b)$
to be
$t_{j}(a,b)$
itself.
Now,
consider
$$
L=\rho_{0}(A)\oplus
\rho_{1}(A)\oplus
\rho_{2}(A)\oplus
T(A,A)\eqno(3.3)$$
where
$T(A,A)$
is a vector space spanned
by $T_{j}(a,b)$
for
any $j=0,1,2$
and for any
$a,b\in A$.
We first assume the
commutation
relation of
$$
[T_{l}(a,b),T_{m}(c,d)]=
-[T_{m}(c,d),T_{l}(a,b)]$$
$$
=T_{m}(t_{l-m}(a,b)c,d)+
T_{m}(c,t_{l-m}(a,b)d))\eqno(3.4)$$
for any
$l,m=0,1,2$
and for
any
$a,b,c,d\in A.$
Then,
it is easy to verify that
it defines a Lie algebra in view of
 Eq.(1.11b).
Of course,
Eq.(3.4)
is also automatically
satisfied
if we identify
$T_{l}(a,b)=
t_{l}(a,b).$
\par
In order to enlarge this Lie algebra,
let us assume
$(i,j,k)$
to be any cyclic permutation
of $(0,1,2)$, 
and let $\gamma_{j}\ (j=0,1,2)$
to be any non-zero constants .
We now assume
\par
(i)
$$
[\rho_{i}(a),\rho_{i}(b)]=
\gamma_{j}\gamma_{k}^{-1}
T_{3-i}(a,b)\eqno(3.5a)$$
\par
(ii)
$$
[\rho_{i}(a),\rho_{j}(b)]=
-[\rho_{j}(b),\rho_{i}(a)]=
-\gamma_{j}\gamma_{i}^{-1}\rho_{k}(\overline{ ab})\eqno(3.5b)$$
\par
(iii)
$$
[T_{l}(a,b),\rho_{j}(c)]=
-[\rho_{j}(c),
T_{l}(a,b)]=
\rho_{j}(t_{l+j}(a,b)c).\eqno(3.5c)$$
\par
Assuming
$\rho_{j}(a)$
to be
$F$-linear in
$a\in A,$
then
$$
L_{j}=T_{3-j}(A,A)
\oplus \rho_{j}(A)\eqno(3.6)$$
yields Lie algebras for each
value of
$j=0,1,2.$
This generalizes 
Eqs.(3.1)
and
(3.2).
Here,
the indices $l$ and $j$ for
$T_{l}(a,b)$
and for $\rho_{j}(a)$
are defined modulo
$3$.
\par
Introducing the Jacobian in $L$
by
$$
J(X,Y,Z)=
[[X,Y],Z]+
[[Y,Z,],X]+
[[Z,X],Y]\eqno(3.7)$$
for
$X,Y,Z\in L,$
we can show
(see Theorem 3.1. of [O.1])
\par
{\bf Proposition 3.1}
\par
\it
Let $A$ be a pre-structurable algebra.
Then,
the Jacobian
$J(X,Y,Z)$
in $L$ are identically zero
except for 
the case of
$$
J(a,b,c):=
J(\rho_{0}(a),\rho_{1}(b),\rho_{2}(c))\eqno(3.8a)$$
which is given by
$$
J(a,b,c)=
T_{0}(a,\overline{ bc})+
T_{1}(c,\overline{ ab})+
T_{2}(b,\overline{ ca}).\eqno(3.8b)$$
\rm
We next note that in view of Eq.(1.11b)
and Theorem 2.1,
we have
$$
[J(a,b,c),\rho_{j}(d)]=
\rho_{j}(Q(a,b,c)d),\eqno(3.9a)$$
$$
[J(a,b,c),T_{j}(d,f)]=
T_{j}(Q(a,b,c)d,f)+
T_{j}(d,Q(a,b,c)f).\eqno(3.9b)$$
\par
Suppose now that $A$ is structurable.
We then have
$Q(a,b,c)=0$
identically,
and
Eqs.(3.9)
yields
$$
[J(a,b,c),\rho_{j}(d)]=
0=
[J(a,b,c),T_{i}(d,f)]\eqno(3.10)$$
so that $J(a,b,c)$
is a center element of $L.$
Let $J$ be a vector space spanned by all
$J(a,b,c),\ (a,b,c\in A).$
Then,
the quotient algebra
${\tilde L}=L/J$
is a Lie algebra.
Therefore,
we can effectively set
$$
J(a,b,c)=T_{0}(a,\overline{ bc})+
T_{1}(c,\overline{ ab})+
T_{2}(b,\overline{ ca})=0.\eqno(3.11)$$
Then,
it is more economical 
to identify
$T_{j}(a,b)$
with
$t_{j}(a,b)$
itself or with a triple of
$$
T(t_{j}(a,b),t_{j+1}(a,b),
t_{j+2}(a,b))\eqno(3.12)$$
as in [A-F],
[E.1]
and [E.2].
In that case,
Eqs.(3.4)
and (3.11) are automatically satisfied by
Eqs.(1.11b)
and Eq.(X).
See also the 
construction given in[K-O].
\par
In what follows,
we suppose now
that $A$ is structurable so that
$L$ is a Lie algebra assuming the valitity of Eq. (3.11).
A special case of
$\gamma_{0}=\gamma_{1}=\gamma_{2}=1$
in Eq.(3.5)
is of particular interest then,
since
$L$ is invariant under a
cyclic permutation group
$Z_{3}$ given by
$$
\rho_{0}(a)\rightarrow
\rho_{1}(a)\rightarrow
\rho_{2}(a)\rightarrow
\rho_{0}(a),\eqno(3.13)$$
$$
T_{0}(a,b)\rightarrow
T_{2}(a,b)\rightarrow
T_{1}(a,b)\rightarrow
T_{0}(a,b).$$
Actually,
$L$ is known to be invariant
under a larger
symmetric group
$S_{4}$
(see [E-O,2], [K-O]),
although we will 
${\underline {\rm not}}$ go into its detail.
\par
We can visualize the structure of the Lie algebra
$L$ given by
Eqs.(3.5)
as in Fig.1, which exhibits a triality:
\par
\vskip 10mm
\begin{center}
\unitlength 0.1in
\begin{picture}(28.06,22.84)(30.60,-35.00)
%
\special{pn 8}%
\special{ar 5290 2200 1100 1100  3.1325839 3.1415927}%
%
\special{pn 8}%
\special{ar 4190 2200 1110 1110  5.2450326 6.2831853}%
\special{ar 4190 2200 1110 1110  0.0000000 2.1332943}%
%
\special{pn 8}%
\special{ar 5300 2190 1084 1084  1.0209431 1.0319768}%
%
\special{pn 8}%
\special{ar 5300 2200 1150 1150  1.0560122 4.1689114}%
%
\special{pn 8}%
\special{ar 4710 3180 1140 1140  3.1769761 6.2831853}%
\special{ar 4710 3180 1140 1140  0.0000000 0.0175421}%
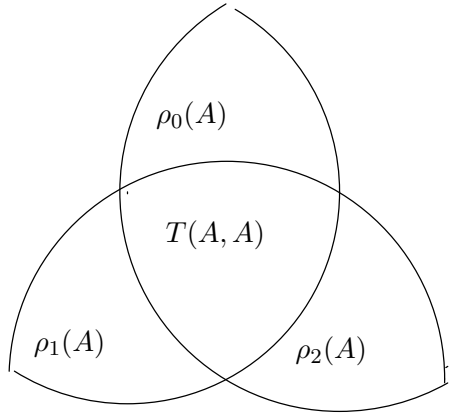
\put(30.6000,-36.7000){\makebox(0,0)[lb]{Fig.1\ {\it Graphical Representation}\ of the Lie Algebra $L$.}}%
\put(43.4000,-18.8000){\makebox(0,0)[lb]{$\rho_{0}(A)$}}%
\put(43.9000,-25.0000){\makebox(0,0)[lb]{$T(A,A)$}}%
\put(37.0000,-30.7000){\makebox(0,0)[lb]{$\rho_{1}(A)$}}%
\put(50.7000,-31.1000){\makebox(0,0)[lb]{$\rho_{2}(A)$}}%
\end{picture}%

\end{center}
\vskip 10mm
\par
In each branch of
the tri-foglio
in Fig.1,
$$
\begin{array}{l}
L_{0}=
\rho_{0}(A)\oplus
T_{0}(A,A)\\
L_{1}=\rho_{1}(A)\oplus
T_{2}(A,A)\\
L_{2}=
\rho_{2}(A)\oplus T_{1}(A,A)
\end{array}
$$
yields these sub-Lie algebra of $L$
as in Eq.(3.6).
They are isomorphic to each other and
interchanges under
the
$Z_{3}$-group
Eq.(3.13)
as in
$$
L_{0}\rightarrow
L_{1}\rightarrow
L_{2}\rightarrow
L_{0}.$$
Also,
$T(A,A)$ is a sub-Lie algebra of $L$,
which transforms among themselves under
$Z_{3}$ as in Eq.(3.13).
\par
We will give some  examples below, assuming underlying
field $F$ to be algebraically closed 
of characterictic zero.
\par
First
let $A$  be an octonion algebra
which is structurable
([A-F]).
In that case,
these Lie algebras are
$L=F_{4},\ L_{j}=B_{4}$
for $j=0,1,2$
and $T(A,A)=D_{4},$
corresponding to the case of  the classical triality relation,
as we may see from works of
([B-S], [E.1] and [E.2]).
\par
If we choose $A$ to be the Zorn's
vector matrix algebra with
Eq.(1.20)
and
(1.21),
corresponding to $B$ beging the 
cubic admissible algebra associated
with the $27$-dimentional exceptional
Jordan algebra
(i.e.Albert algebra),
then $A$ is structurable
and the resulting
Lie algebra $L$ is of type $E_{8}$
(see e.g. [Kan] and [Kam.1]).
Moreover,
$L_{j}(j=0,1,2)$
is a Lie algebra
$E_{7}\oplus A_{1}$
and
$T(A,A)$
is realized
to be isomorhic to
$E_{6}\oplus gl(1)\oplus gl(1).$
Another way of obtaining
$E_{8}$ is to consider a structurable algebra
$A=O_{1}\otimes O_{2}$
of two octonion algebra
$O_{1}$
and $O_{2}$
(see
[A-F]).
In that case,
it is known
(see
[B-S],[E.1]
and [E.2])
that it yields also the
Lie algebra
$E_{8}$.
However,
the sub-Lie algebra
$L_{j}(j=0,1,2)$
in this case are
Lie algebra
$D_{8},$
while
$T(A,A)$ is
$D_{4}\oplus D_{4}.$
\par
It may be worth-while to make the following comment here.
The Lie algebra constructed by Eqs.(3.3)-(3.6) manifests the 
explicit $Z_{3}$ -symmetry (i.e. the triality ), but $\underline{not}$
 a 5-graded structure. On the other side,
the standard construction of the Lie algebra on the basis of
the (-1,1) Freudenthal-Kantor triple system[Y-O] is on the contrary explicitly
5-graded but $\underline{not} $ manifestly $Z_{3}$-invariant
.  
The relationship between these two approaches has been studied in [K-O].

\vskip 3mm

In ending of this note, we remark
\par
{\bf Remark 3.2}
\par
\vskip 1mm
Any simple structurable algebra
$A$
may be identified with
some symmetric space as follows.
First,
if $A$ is structurable,
then we can construct 
Lie algebras
by Eqs.(3.5)
with Eqs.(3.4) and
(3.13)
from the standard construction based upon
$(-1,1)$Freudenthal-Kantor triple system
[K-O].Then by
Eqs.(3.5),
we have
$$
[\rho_{0}(a),
\rho_{0}(b)]
=
(\gamma_{2}/\gamma_{0})T_{0}(a,b)
$$
$$
[T_{0}(a,b),\rho_{0}(c)]
=
\rho_{0}(t_{0}(a,b)c)
$$
$$
[T_{0}(a,b),T_{0}(c,d)]
=
T_{0}(t_{0}(a,b)c,d)+
T_{0}(c,t_{0}(a,b)d)
$$
so that
$\rho_{0}(A)$
may be identified with the symmetric
space
$L_{0}/T_{0}(A,A).$
because 
 the tangent space of a symmetric space 
has a structure of Lie triple system.
Further,
the mapping
$A\rightarrow \rho_{0}(A)$
is one-to-one
if $A$ is simple.
To prove it,
let
$$
B=\{a|\rho_{0}(a)=0,\ a\in A\}
$$
and
calculate
$$
\begin{array}{l}
0=[\rho_{0}(a),\rho_{1}(x)]=
-(\gamma_{1}/\gamma_{0})
\rho_{2}(\overline{ ax})\\
0=
[\rho_{2}(x),\rho_{0}(a)]=
-(\gamma_{0}/\gamma_{2})\rho_{1}
(\overline{ xa})
\end{array}
$$
for any
$a\in B$
and any
$x\in A$,
so that
$\rho_{1}(\overline{ xa})=
\rho_{2}(\overline{ ax})=0.$
Moreover,
we compute
$$
\begin{array}{l}
0=[\rho_{1}(e),\rho_{2}(\overline{ ax})]=
-(\gamma_{2}/\gamma_{1})
\rho_{0}(ax)\\
0=[\rho_{1}(\overline{ xa}),
\rho_{2}(e)]=
-(\gamma_{1}/\gamma_{2})
\rho_{0}(xa)
\end{array}
$$
which yields
$\rho_{0}(ax)=\rho_{0}(xa)=0.$
Therefore
$B$ is a ideal of $A$ and hence
$A\rightarrow \rho_{0}(A)$
is one-to-one,
if $A$ is simple with
$B=0.$
Note that the other possibility of $B=A$ leads to the trivial
case of $\rho _{0}(A)=T_{0}(A,A) =0$ identically.

\vskip 3mm
{\bf References}
\par
\rm
\vskip 3mm
\noindent
[A]:B.N.Allison;
A class of non-associative algebra with involution
containing
the class of Jordan algebra,
Math,Ann 
${\bf {\rm 237}}$,
(1978)
133-156
\par
\noindent
[A-F]:
B.N.Allison
and J.R.Faulkner;
Non-associative algebras for
Steinberg
unitary Lie algebras,
J.Algebras ${\bf {\rm 161}}$
(1993)1-19
\par
\noindent
[B-S]:
C.M.Burton and A.Sudbury;
Magic Squares
and matrix module of Li algebras,
Adv.Math.
${\bf {\rm 180}}$
(2003)
596-647
\par
\noindent
[E.1]:A.Elduque;
The magic square and symmetric compositions,
Rev.
Math.
Iberoamerica
${\bf {\rm 20}}$(2004)
477-493
\par
\noindent
[E.2]: A.Elduque,
A New look of Freudenthal- magic square.
In {\it Non-associative Algebras and 
 its applications},
ed.by L.Sabinin,
L.Sbitneva,
and I.Shestakov,
Chap and Hall,
new York
(2006),
pp149-165
\par
\noindent
[E-O.1]:
A.Elduque, and S.Okubo;
On algebras satisfying
$x^{2}x^{2}=N(x)x,$
Math.Zeit.
${\bf {\rm 235}}$
(2000)
275-314
\par
\noindent
[E-O.2]
A.Elduque,
and S.Okubo;
Lie algebras with $S_{4}$-action and
structurable algebras,
J.Algebras
${\bf {\rm 307}}$
(2007)
864-890
\par
\noindent
[E-K-O];
A.Elduque,N.Kamiya,
and S.Okubo;
Left unital Kantor triple system and
structurable algebra,
arXiv
1205,
2489(2012), to appear in Linear and Multi-linear algebras.
\par
\noindent
[Kam.1]:
N.Kamiya;
A structure theory of Freudenthal-Kantor
triple 
system 
III,
Mem.
Fac.
Sci.
Shimane 
Univ.
${\bf{\rm 23}}$
(1989)
33-51
\par
\noindent
[Kam.2]:
N.Kamiya;
On generalized structurable algebras and
Lie-related 
triple,
Adv.
Clifford Algebras
${\bf {\rm 5}}$
(1995)
127-140
\par
\noindent
[Kan]:
I.L.Kantor,
Models of the exceptional
 Lie algebras,
 Sov.
Math.
Dokl.
${\bf {\rm 14}}$(1973)
254-258
\par
\noindent
[Ku]:
E.N.Kuzumin;
Malcev-algebras of dimension
five over
a field of Zero characteristic,
Algebra
i Logika
${\bf{\rm 9}}$
(1970)
691-700
$=$
Algebra and Logic 
${\bf {\rm 9}}$
(1970)
416-421
\par
\noindent
[K-O]:
N.Kamiya and S.Okubo, 
Symmetry of Lie algebras associated with 
$(\varepsilon,\delta)$-Freudenthal-Kantor 
triple systems, arXiv 1303,0072.
\par
\noindent
[M]:
H.C.Myung;
{\it Malcev-Admisible\ 
Algebras},
\rm
Birkhauser,
Boston
(1986)
\par
\noindent
[O.1]:
S.Okubo;
Symmetric triality relations and 
structurable algebra
Linear Algebras and
its Applications,
{\bf {\rm 396}}
(2005)
189-222
\par
\noindent
[O.2]:
S.Okubo;
{\it Algebras satisfying symmetric 
triality relations,}
in
 {\rm Non-associative
\ Algebras}
\par
\noindent
{\bf{\rm and\  
its\ 
Applications,}}
ed.by L.Sabinin,
L.Sbitneva,
and 
I.P.Shestakov,
Chapman 
and Hall.
N.Y.
(2006)
313-321
\par
\noindent
[P-S]:
J.M.Perez-Izquierdo and
I.P.Shestakov;
An envelope for Malcev algebra
J.Algebra,
{\bf {\rm 272}}
(2004)
379-393
\par
\noindent
[S]:
R.D.Schafer;
${\it An\ Introduction\ 
to \ Non-associative\ algebras}$,
Academic Press,
New York
(1966)
\par
\noindent
[Y-O]:
K.Yamaguti
and A.Ono;
On representation of Freudenthal-Kantor 
triple system
$U(\varepsilon,\delta)$,
Bull.Fac.School.
Ed.Hiroshima Univ.
Part II,
${\bf {\rm 7}}$
(1984)43-51
\end{document}